\documentclass{amsart}

\usepackage{amsthm}
\usepackage{graphicx}

\newtheorem{theo}{Theorem}

\def\cB{\mathcal{B}}
\def\cE{\mathcal{E}}
\def\cT{\mathcal{T}}
\def\cBs{\cB_{\sigma}}
\def\cEs{\cE_{\sigma}}
\def\cEsp{\cE_{\sigma'}}
\def\ni{\noindent}
\def\ms{\medskip}
\def\bn{\mathbf{n}}
\def\Snij{S_n^{(i,j)}}

\title{On symmetries in phylogenetic trees}
\author[\'E. Fusy]{\'{E}ric Fusy$^{*}$}
\thanks{$^{*}$LIX, \'Ecole Polytechnique, Palaiseau, France, fusy@lix.polytechnique.fr. Partly supported by the ANR grant ``Cartaplus'' 12-JS02-001-01 and the ANR grant ``EGOS'' 12-JS02-002-01.}

\begin{document}

\begin{abstract}
Billey et al. [arXiv:1507.04976] have recently discovered a surprisingly simple
formula for the number $a_n(\sigma)$ 
of leaf-labelled rooted non-embedded binary trees
(also known as phylogenetic trees) with $n\geq 1$ leaves,   
fixed (for the relabelling action)
by a given permutation $\sigma\in\frak{S}_n$. 
Denoting by $\lambda\vdash n$ the integer partition giving the 
sizes of the cycles of $\sigma$ in non-increasing order, they show by a 
guessing/checking approach that if $\lambda$ is a binary partition (it is known that $a_n(\sigma)=0$ otherwise), then
$$
a_n(\sigma)=\prod_{i=2}^{\ell(\lambda)}(2(\lambda_i+\cdots+\lambda_{\ell(\lambda)})-1),
$$
and they derive from it a formula and random generation procedure for tanglegrams
(and more generally for tangled chains).  Our main result is a combinatorial proof
of the formula for $a_n(\sigma)$, 
which yields a simplification of the random sampler for 
tangled chains.
\end{abstract}

\maketitle

\section{Introduction}
\begin{figure}[h!]
\begin{center}
\includegraphics[width=9cm]{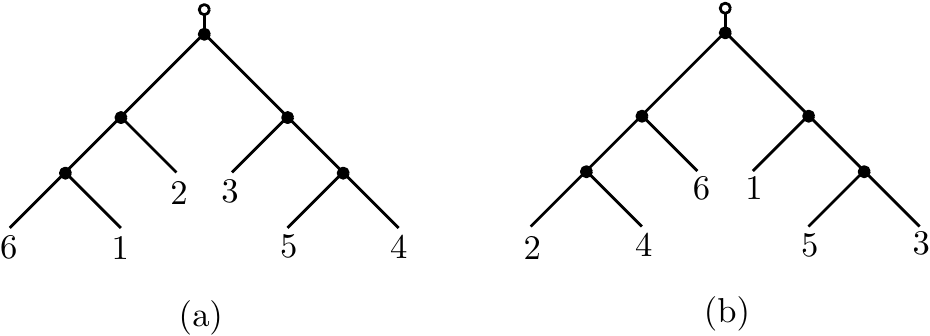}
\end{center}
\caption{(a) A phylogenetic tree $\gamma$ with label-set $[1..6]$. 
(b) The tree $\gamma'=\sigma\cdot\gamma$, with $\sigma=(1,4,3)(5)(2,6)$.
Since $\gamma'\neq \gamma$, $\gamma$ is not fixed by $\sigma$ (on the
other hand $\gamma$ is fixed by $(2,3)(1,4,6,5)$).}
\label{fig:examples}
\end{figure}

For $A$ a finite set of cardinality $n\geq 1$, we denote by $\cB[A]$ the set of 
rooted binary trees that are non-embedded (i.e., the order of 
the two children of each node does not matter) and have $n$ leaves with distinct
labels from $A$. Such trees are known as \emph{phylogenetic trees}, 
where typically $A$ is the set of represented species.  
Note that such a tree has $n-1$ nodes and $2n-1$ edges
(we take here the convention of having an additional root-edge above the root-node, connected
to a `fake-vertex' that does not count as a node, 
see Figure~\ref{fig:examples}). 

The group $\frak{S}(A)$ of permutations of $A$ acts on $\cB[A]$:
for $\gamma\in\cB[A]$ and $\sigma\in\frak{S}(A)$, 
$\sigma\cdot \gamma$ is obtained from $\gamma$ after replacing
 the label $i$ of every leaf by $\sigma(i)$, see Figure~\ref{fig:examples}(b).  
We denote by 
$\cBs[A]$ 
  the set of trees fixed by the action of $\sigma$, i.e., $\cBs[A]:=\{\gamma\in\cB[A]\ \mathrm{such\ that}\ \sigma\cdot\gamma=\gamma\}$.
 We also define $\cEs[A]$ (resp. $\cE[A]$) 
as the set of pairs $(\gamma,e)$
where $\gamma\in\cBs[A]$ (resp. $\gamma\in\cB[A]$)
 and $e$ is an edge of $\gamma$ (among the $2n-1$ edges).  
Define the \emph{cycle-type} of $\sigma$ as the integer
partition $\lambda\vdash n$ giving the sizes of the cycles of $\sigma$
(in non-increasing order). 
For $\lambda\vdash n$ an integer partition, the cardinality of $\cBs[A]$
is the same for all permutations $\sigma$ with cycle-type $\lambda$,
and this common cardinality is denoted by $r_{\lambda}$.  
It is known (e.g. using cycle index sums~\cite{BeLaLe,gessel2015}) 
that $r_{\lambda}=0$ 
 unless $\lambda$ is a binary partition (i.e., an integer partition whose parts are powers of $2$).  
Billey et al.~\cite{billey2015} have recently found  
the following remarkable formula, valid  
for any binary partition $\lambda$:
\begin{equation}\label{eq:rlambda}
r_{\lambda}=\prod_{i=2}^{\ell(\lambda)}(2(\lambda_i+\cdots+\lambda_{\ell(\lambda)})-1).
\end{equation} 
They prove the formula by a guessing/checking approach. 
Our main result here is a combinatorial proof of~\eqref{eq:rlambda}, which   
 yields a simplification (see Section~\ref{sec:random_gen}) 
of the random sampler for tanglegrams (and more
generally tangled chains) given in~\cite{billey2015}.

\begin{theo}\label{theo:m}
For $A$ a finite set and $\sigma$ a permutation on $A$ whose cycle-type is a binary partition:
\begin{itemize}
\item
If $\sigma$ has one cycle, then $|\cBs[A]|=1$.
\item
If $\sigma$ has more than one cycle, let $c$ be a largest 
cycle of $\sigma$; denote by $A'$ the set $A$ without the 
elements of $c$, and denote by $\sigma'$ the permutation
$\sigma$ restricted to $A'$. Then we have the combinatorial isomorphism
\begin{equation}\label{eq:identity}
\cBs[A]\simeq \cEsp[A'].
\end{equation}
\end{itemize}
\end{theo}

As we will see, the isomorphism~\eqref{eq:identity} 
can be seen as an adaptation of R\'emy's method~\cite{remy1985procede}    
to the setting of (non-embedded rooted) binary trees fixed by a given permutation.   
Note that Theorem~\ref{theo:m} implies that the coefficients $r_{\lambda}$ satisfy
 $r_{\lambda}=1$ if $\lambda$ is a binary partition with one part
and $r_{\lambda}=(2|\lambda\backslash \lambda_1|-1)\cdot r_{\lambda\backslash \lambda_1}$ if
$\lambda$ is a binary partition with more than one part, from
which we recover~\eqref{eq:rlambda}.

\section{Proof of Theorem~\ref{theo:m}} 
\subsection{Case where the permutation $\sigma$ has one cycle}\label{sec:one}
The fact that $|\cBs[A]|=1$ if $\sigma$ has one cycle of size $2^k$ (for some $k\geq 0$)  
is well known from the structure of automorphisms in trees~\cite{polya1937kombinatorische},   
for the sake of completeness we give a short justification.  
Since the case $k=0$ is trivial we can assume that $k\geq 1$.
 Let $c_1,c_2$ be the two cycles of $\sigma^2$ (each of size $2^{k-1}$), with the convention that $c_1$ 
contains the minimal element of $A$; denote by $A_1,A_2$ the induced bi-partition of $A$,
and by $\sigma_1=c_1$ (resp. $\sigma_2=c_2$) the permutation $\sigma^2$ restricted to $A_1$ (resp. $A_2$).  
  For $\gamma\in\cBs[A]$ let $\gamma_1,\gamma_2$ be the two subtrees at the root-node of $\gamma$,
such that the minimal element of $A$ is in $\gamma_1$. Then clearly $\gamma_1\in\cB_{\sigma_1}[A_1]$
and  $\gamma_2\in\cB_{\sigma_2}[A_2]$, and conversely for $\gamma_1\in\cB_{\sigma_1}[A_1]$
and  $\gamma_2\in\cB_{\sigma_2}[A_2]$ the tree $\gamma$ with $(\gamma_1,\gamma_2)$ as subtrees at the root-node
is in $\cBs[A]$. Hence
\begin{equation}
\cBs[A]\simeq \cB_{\sigma_1}[A_1]\times\cB_{\sigma_2}[A_2],
\end{equation}
which implies $|\cBs[A]|=1$ by induction on $k$ (note that, also 
by induction on $k$, the underlying unlabelled tree is the complete binary tree of height $k$).

\subsection{Case where the permutation $\sigma$ has more than one cycle} 
Let $k\geq 0$ be the integer such that the largest cycle of $\sigma$ has size $2^k$. 
A first useful remark is that $\sigma$ induces a permutation of the edges (resp. of the nodes) of $\gamma$, and each $\sigma$-cycle
of edges (resp. of nodes) has size $2^i$ for some $i\in [0..k]$. We present the proof of~\eqref{eq:identity} progressively, treating first 
the case $k=0$, then $k=1$, then general $k$. 

\begin{figure}
\begin{center}
\includegraphics[width=12cm]{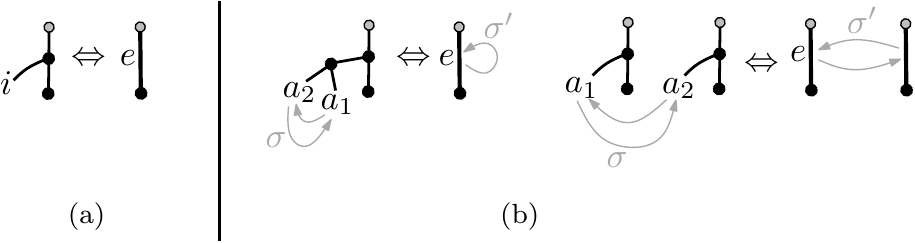}
\end{center}
\caption{(a) R\'emy's leaf-removal operation. (b) The two cases
for removing a $2$-cycle of leaves (depending whether the 
two leaves have the same parent or not). The vertices depicted gray 
are allowed to be the fake vertex above the root-node.}
\label{fig:remy}
\end{figure}

\ms

\ni {\bf Case $k=0$.} This case corresponds to $\sigma$ being the identity, so that $\cBs[A]\simeq\cB[A]$, hence
we just have to justify that $\cB[A]\simeq\cE[A\backslash\{i\}]$ for each fixed $i\in A$. This is easy
to see using R\'emy's argument~\cite{remy1985procede}~\footnote{A similar argument in the context of triangulations of a polygon dates back to Rodrigues~\cite{rodrigues1838nombre}.}, 
used here in the non-embedded leaf-labelled context: every $\gamma\in\cB[A]$ is uniquely obtained from some
 $(\gamma',e)\in\cE[A\backslash\{i\}]$ upon inserting a new pending
edge from the middle of $e$ to a new leaf that is given label $i$, 
see Figure~\ref{fig:remy}(a). 

\ms

\ni {\bf Case $k=1$.} 
Let $c=(a_1,a_2)$ be the selected cycle of $\sigma$, with $a_1<a_2$. 
Two cases can arise (in each case we obtain from
$\gamma$ a pair $(\gamma',e)$ with $\gamma'\in\cB_{\sigma'}[A']$ 
 and $e$ an edge of $\gamma'$): 
\begin{itemize}
\item
if $a_1$ and $a_2$ have the same parent $v$, we obtain a reduced
tree $\gamma'\in\cB_{\sigma'}[A']$ by erasing  
the $3$ edges incident to $v$ (and the endpoints of these edges,
which are $a_1,a_2,v$ and the parent of $v$),
and we mark the edge $e$ of $\gamma'$ whose middle was the parent of $v$,
see the first case of Figure~\ref{fig:remy}(b) 
\item
if $a_1$ and $a_2$ have distinct parents, we can apply the operation 
of Figure~\ref{fig:remy}(a) to each of $a_1$ and $a_2$, which yields a 
reduced tree $\gamma'\in\cB_{\sigma'}[A']$. 
We then mark the edge $e$ of $\gamma'$ 
whose middle was the parent of $a_1$, 
see the second case of Figure~\ref{fig:remy}(b).
\end{itemize}

Conversely, starting from $(\gamma',e)\in\cE[A']$,  
the $\sigma'$-cycle of edges that contains $e$ has either
size $1$ or $2$:
\begin{itemize}
\item
if it has size $1$ (i.e., $e$ is fixed by $\sigma'$), we
insert a pending edge from the middle of $e$ and leading
to  ``cherry" with labels $(a_1,a_2)$,
\item
if it has size $2$, let $e'=\sigma'(e)$; then we attach at the 
middle of  $e$ (resp. $e'$) a new pending edge leading to 
a new leaf of label $a_1$ (resp. $a_2$). 
\end{itemize}

\ms

\ni {\bf The general case $k\geq 0$.} 
\begin{figure}
\begin{center}
\includegraphics[width=12cm]{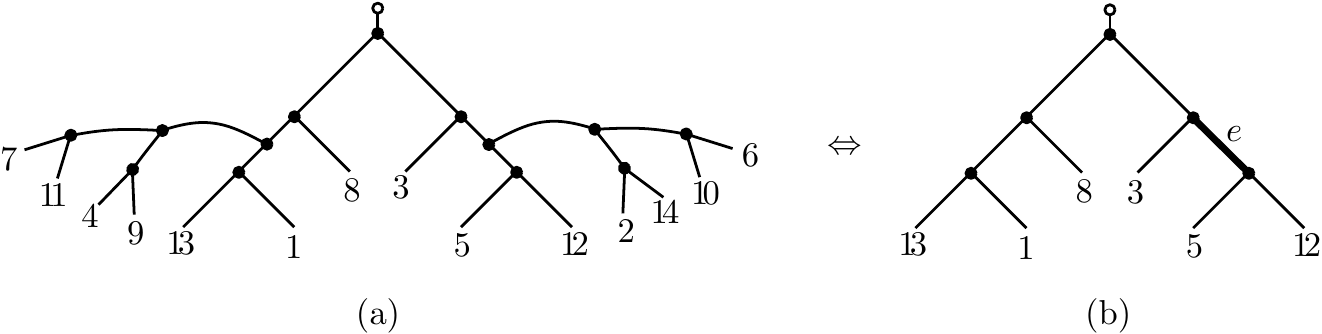}
\end{center}
\caption{(a) a tree in $\cBs[A]$, for $A=[1..14]$ and
$\sigma=(3,8)(1,5,13,12)(2,7,10,4,14,11,6,9)$. 
(b) The corresponding (when selecting the cycle $c$ of size $8$ in $\sigma$) pair $(\gamma',e)\in\cE_{\sigma'}[A']$, where $A'=A\backslash c$ 
and $\sigma'=(3,8)(1,5,13,12)$ (restriction of $\sigma$ to $A'$).}
\label{fig:general}
\end{figure}
Recall that the marked cycle of $\sigma$ is denoted by $c$. 
A node or leaf of the tree is generically called a \emph{vertex}
of the tree.  
We define a \emph{$c$-vertex} as a vertex $v$ of $\gamma$ such 
that:
\begin{itemize}
\item
if $v$ is a leaf then $v\in c$,
\item
if $v$ is a node then 
 all leaves that are descendant of $v$ are in $c$.
\end{itemize}    
A $c$-vertex is called \emph{maximal} if it is not the descendant
of any other $c$-vertex; define a \emph{$c$-tree} as a subtree
formed by a maximal $c$-vertex $v$ and its hanging subtree
(if $v$ is a leaf then the corresponding $c$-tree is reduced to $v$).  
 Note that the maximal $c$-vertices are permuted by $\sigma$.  
Moreover since the leaves of $c$ are permuted cyclically, the
maximal $c$-vertices actually have to form a $\sigma$-cycle 
of vertices,  
of size $2^i$ for some $i\leq k$; and in each $c$-tree,  
$\sigma^{2^i}$ 
permutes the $2^{k-i}$ leaves of the $c$-tree cyclically.  
 Let $\ell$ be the leaf of minimal label in $c$, and let
$w$ be the maximal $c$-vertex such that the $c$-tree at $w$ 
 contains $\ell$.  
 We obtain a reduced tree 
$\gamma'\in\cB_{\sigma'}[A']$ by erasing all $c$-trees and erasing 
the parent-edges
and parent-vertices of all maximal $c$-vertices; and then
we mark the edge $e$ of $\gamma'$ whose middle 
was the parent of $w$, see Figure~\ref{fig:general}. 

Conversely, starting from $(\gamma',e)\in\cE_{\sigma'}[A']$, 
let $i\in[0..k]$ be such that the $\sigma'$-cycle of edges that contains
$e$ has cardinality $2^i$; write this cycle as 
$e_0,\ldots,e_{2^i-1}$, with $e_0=e$. 
Starting from the element of $c$ of minimal label, 
let $(s_0,\ldots,s_{2^i-1})$ be the $2^i$ (successive) first
 elements of $c$.  
And for $r\in[0..2^i-1]$ 
let $c_r$ be the cycle of $\sigma^{2^i}$ that contains $s_r$,
and let $A_r$ be the set of elements in $c_r$  
(note that $A_0,\ldots,A_{2^i-1}$ each have size $2^{k-i}$ and partition
the set of elements in $c$). 
Let $T_r$ be the unique (by Section~\ref{sec:one}) tree in $\cB[A_r]$
fixed by the cyclic permutation $c_r$. We obtain a tree $\gamma\in\cBs[A]$
as follows: for each $r\in[0..2^i-1]$ we create a new edge that connects
the middle of $e_r$ to a new copy of $T_r$.

To conclude we have described a mapping 
from $\cBs[A]$ to $\cEsp[A']$ and a mapping from $\cEsp[A']$ 
 to $\cBs[A]$ that are readily seen to be 
 inverse of each other, therefore  
$\cBs[A]\simeq \cEsp[A']$. 

\section{Application to the random generation of tangled chains}\label{sec:random_gen}

\begin{figure}
\begin{center}
\includegraphics[width=12cm]{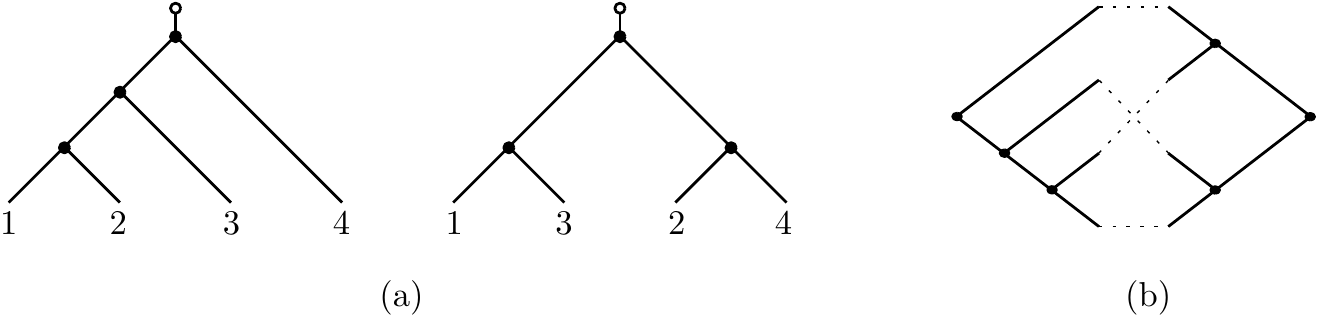}
\end{center}
\caption{(a) A pair of (rooted non-embedded leaf-labelled) binary trees. (b)
the corresponding (unlabelled) tanglegram.}
\label{fig:tangle}
\end{figure}

For $n\geq 1$, denote by $\bn$ the set $\{1,\ldots,n\}$. 
A \emph{tanglegram} of size $n$ is an orbit of $\cB[\bn]\times\cB[\bn]$
under the relabelling action of $\frak{S}_n$ (see Figure~\ref{fig:tangle} for 
an example). More generally, for $k\geq 1$, a \emph{tangled chain} 
of length $k$ and size $n$ is an orbit of $\cB[\bn]^k$
under the relabelling action of $\frak{S}_n$, see~\cite{matsen2015tanglegrams, billey2015, gessel2015}.  
Let $\cT_{n}^{(k)}$ be the set 
of tangled chains of length $k$ and size $n$, and let
$t_{n}^{(k)}$ be the cardinality of $\cT_{n}^{(k)}$. 
Then it follows from Burnside's lemma (see~\cite{billey2015} for a proof using 
double cosets and~\cite{gessel2015} for a proof using the formalism of species) 
 that
\begin{equation}\label{eq:tnk}
t_{n}^{(k)}=\frac1{n!}\sum_{\sigma\in\frak{S}_n}|\cBs[\bn]|^k=\sum_{\lambda\vdash n}\frac{r_{\lambda}\ \!\!\!^k}{z_{\lambda}}, 
\end{equation}
where $z_{\lambda}=1^{m_1}m_1!\cdots r^{m_r}m_r!$ if $\lambda$ has $m_1$ parts
of size $1$,...,$m_r$ parts of size $r$ (recall that $n!/z_{\lambda}$ is the 
number of permutations with cycle-type $\lambda$).  
At the level of combinatorial classes, Burnside's lemma gives
$$
\frak{S}_n\times\cT_{n}^{(k)}\simeq \sum_{\sigma\in\frak{S}_n}\cBs[\bn]^k,
$$
and thus the following procedure is a uniform random sampler for $\cT_n^{(k)}$ (see~\cite{billey2015} for details):
\begin{enumerate}
\item
Choose a random binary partition $\lambda\vdash n$ under the distribution
$$
P(\lambda)=\frac{r_{\lambda}\ \!\!\!^k/z_{\lambda}}{S_n},
$$
where $S_n=\sum_{\lambda\vdash n}r_{\lambda}\ \!\!\!^k/z_{\lambda}$ ($=t_n^{(k)}$).
\item
Let $\sigma$ be a permutation with cycle-type $\lambda$. 
For each $r\in[1..k]$ draw (independently) a tree $T_r\in\cB_{\sigma}[\bn]$
uniformly at random. 
\item
Return the tangled chain corresponding to $(T_1,\ldots,T_k)$. 
\end{enumerate}  
A recursive procedure (using~\eqref{eq:rlambda})  
is given in~\cite{billey2015} to sample
uniformly at random from $\cB_{\sigma}[\bn]$. From Theorem~\ref{theo:m} we obtain 
a simpler random sampler for $\cB_{\sigma}[\bn]$. 
We order the cycles of $\sigma$ as  
$c_1,\ldots,c_{\ell(\lambda)}$ such that the cycle-sizes are in 
non-decreasing order. 
Then, with $A_1$ the set of labels in $c_1$, 
 we start from the unique tree (by Section~\ref{sec:one}) in $\cB_{c_1}[A_1]$
(where $c_1$ is to be seen as a cyclic permutation on $A_1$). 
Then, for $i$ from $2$ to $\ell(\lambda)$ we mark an edge chosen 
uniformly at random from the already obtained tree, and then 
we insert the leaves that have labels in $c_i$ 
using the isomorphism~\eqref{eq:identity}.

The complexity of the sampler for $\cB_{\sigma}[\bn]$ is clearly linear in $n$ and needs no precomputation of coefficients. 
However step~(1) of the random generator 
requires a table of $p(n)$ coefficients, where $p(n)$ is the number of 
binary partitions of $n$, 
which is slightly superpolynomial~\cite{mahler1940special}, 
$p(n)=n^{\Theta(\log(n))}$.
It is however possible to do step~(1) in polynomial time.
For this, we consider, for $i\geq 0$ and $n,j\geq 1$ 
the coefficient $\Snij$ defined as the 
sum of $r_{\lambda}\ \!\!\!^k/z_{\lambda}$ over all binary partitions
of $n$ where the largest part is $2^i$ and has multiplicity $j$;
note that $\Snij=0$ unless $j\cdot 2^i\leq n$, we denote by $E_n$
the set of such pairs $(i,j)$. Since $r_{\lambda}=1$ and $z_{\lambda}=(|\lambda|-1)!$  if $\lambda$ has one part,
we have the initial condition $\Snij=1/(n-1)!$ for $j=1$ and $2^i=n$. 
In addition, using the fact that 
$r_{\lambda}=(2|\lambda\backslash \lambda_1|-1)\cdot r_{\lambda\backslash \lambda_1}$ if $\lambda$ has at least $2$ parts, and the formula for $z_{\lambda}$, 
we easily obtain the recurrence:
$$
\Snij=\frac{(2(n-2^i)-1)^k}{2^ij}S_{n-2^i}^{(i,j-1)}\ \mathrm{for}\ (i,j)\in E_n\ \mathrm{with}\ 2^i<n,
$$
valid for $j=1$ upon defining by convention 
 $S_n^{(i,0)}$ as the sum of $S_{n}^{(i',j')}$ over all pairs $(i',j')\in E_{n}$ such that $i'<i$.  

Thus in step~(1), instead of directly drawing $\lambda$ under $P(\lambda)$,  
we may first choose the pair $(i,j)$ such that the largest part of $\lambda$ 
is $2^i$ and has multiplicity $j$, that is, we draw $(i,j)\in E_n$ under distribution $P(i,j)=\Snij/S_n$. 
Then we continue recursively at size $n'=n-2^ij$, but
 conditioned on the largest part to be smaller than $2^i$
(that is, for the second step and similarly for later steps, we draw the  pair $(i',j')$ in $E_{n'}\cap\{i'<i\}$ 
under distribution $S_{n'}^{(i',j')}/S_{n'}^{(i,0)}$).  
Note that 
$|E_n|=\sum_{i\leq \log_2(n)}\lfloor n/2^i\rfloor=\Theta(n)$. Since 
we need all coefficients $S_{m}^{(i,j)}$ for $m\leq n$ and $(i,j)\in E_m$,
we have to store $\Theta(n^2)$ coefficients. In addition it is easy to see (looking at the first expression in~\eqref{eq:tnk}) 
that each coefficient $S_m^{(i,j)}$ is a rational number of the form $a/m!$
with $a$ an integer having $O(m\log(m))$ bits. Hence the overall storage 
bit-complexity 
is $O(n^3\log(n))$. About time complexity, starting at size $n$
we first choose the pair $(i,j)$ (with $2^i$ the largest part and $j$
its multiplicity), which takes  $O(|E_n|)=O(n)$ comparisons, and then 
we continue recursively at size $n-j\cdot 2^i$. 
At each step the choice of a pair $(i,j)$ takes time $O(m)$ with $m\leq n$ the 
current size,  and the number
of steps is the number of distinct part-sizes in the finally 
output binary partition 
$\lambda\vdash n$. Since the number of distinct part-sizes in a binary 
partition of $n$ is $O(\log(n))$, we conclude that the time complexity 
(in terms of the number of real-arithmetic comparisons) 
to draw $\lambda$ is $O(n\log(n))$.

\ms
\ms

\ni\emph{Acknowledgements.} I thank Igor Pak
for interesting discussions.

\bibliographystyle{plain}
\bibliography{mabiblio}

\end{document}